\documentclass{article}
\usepackage{maa-empty}

\usepackage{amsthm}

\theoremstyle{definition}

\newcommand{\old}[1]{}

\begin{document}

\title{What is a proof? What should it be?}
\markright{Abbreviated Article Title}
\author{Christoph Benzm{\"u}ller}

\maketitle

\vspace*{-1em}
\begin{abstract}
\emph{Mathematical proofs} should be paired with \emph{formal proofs}, whenever feasible.
\end{abstract}

% MSC classes: 03B35, 68T15, 97D20, 68N99, 68T15, 68T30
% ACM classes: F.4; I.2.3; I.2.4; D.2.4; G.4; K.3.1; J.0

What is a proof? Is it the rigorous but typically rather unintuitive \emph{formal derivation of a new ``truth'' from its premises using accurately defined rules of inference}?
%, as provided, for example, in Gentzen style natural deduction or sequent calculi
Or is it an \emph{artful communication act} in which the beautiful structures underlying a new mathematical insight are revealed to peer experts in such a way that they can easily \emph{see} and accept it,
%the new ``truth'', 
and even gain further inspiration?

% At first sight, the contrast between these two antipodes could hardly be bigger.

The former notion, referred to as  \emph{formal proof}, is primarily concerned with logical rigor and soundness. Intuition and beauty is still a secondary concern, if at all. Formal proofs have recently attained increased, albeit quite controversial, attention in mathematics. This interest has been triggered by successful applications of modern theorem proving technology to challenging mathematical verification and reasoning tasks. Some of the settled problems are of such a kind, that human cognition alone has apparently reached its limits
%boundaries 
for attacking them. Respective examples include:
\begin{enumerate}
    \item \label{hales} Hales' \cite{Hales2017} verification of his proof of the Kepler conjecture within the proof assistant HOL Light: A board of expert reviewers of the Annals of Mathematics had previously surrendered this complex task, but Hales and his team mastered it in interaction with a proof assistant system. As the main result,  the computer system produced a formal proof that is now independently verifiable --- by humans and/or (other) computer programs.
    \item \label{heule} Heule \& Kullmann's \cite{Heule2017} automated solution of the Pythagorean Triples Problem: In this work an open mathematical problem was solved with fully automated SAT solving technology. The formal proof that was generated by the computer program  is of enormous size (about 200TB). Nevertheless, it is still  independently verifiable (at least by machines). This line of research has recently been continued by an automatic solution for Schur Number Five \cite{Heule2018}.
\end{enumerate}
While some mathematicians embrace this new, computer-supported alternative mathematics, many others still strongly reject it and ask: ``Is this still maths?''. 

Those latter, disapproving mathematicians typically point to the virtues of traditional \emph{mathematical proofs}, which, in contrast to formal proofs, focus on intuition, beauty and explanatory power. However, such proofs often lack logical rigor, and the exact dependencies and the precisely required inference principles may remain vague --- but these weaknesses are considered subordinate to human intuition and abstract-level understanding. Moreover, the assessment of mathematical proofs is traditionally also handled quite differently from those of formal proofs: it is organized as a kind of social/voting process in which a sufficient number of peers has to be convinced of the new result for it to be generally established. 
It is thus of little surprise that many mathematical proofs do in fact suffer from mostly minor, but occasionally also major, technical flaws that have escaped the human eye. Most of these errors, I claim, would have been revealed in a formal proof verification process. 

Both notions of proof thus constitute strong antipodes, with orthogonal pros and cons, and with opposing goals. Traditional mathematical proofs are made for, and consumed by, humans, while formal proofs are predominantly generated with, and consumed by, machines.

\emph{So, what should a proof be?} In my opinion it should ideally be both, \emph{whenever feasible}, namely a human-oriented traditional proof accompanied by a machine-oriented formal proof.  
There might be situations though in which only one of both notions can be provided (in principle or for the time being). For example, it is still unclear whether the 200TB proof generated in \eqref{heule} can be replaced by, respectively accompanied with, a human-oriented, short and intuitive proof --- simply because there might be none. In fact, due to the sprouting complexity in an increasingly technified world, we can actually expect a soaring number of analogous challenges to emerge in the future, in particular, in areas such as computer science and artificial intelligence. Think e.g.~about the assessment and verification of critical software components in emerging intelligent systems. We cannot even expect beautiful and insightful proofs to generally exist in such contexts, since the systems to be assessed might simply be too complex while at the same time ill-designed or relying on ill-defined foundations. But can we, or should we, therefore capitulate from verification attempts, only because human intuitive proofs or refutations are not easily feasible in a particular application context? Clearly not! I am convinced that we even have the duty to take on such challenges. I am thus strongly against upholding a restricted, traditional notion of mathematical proof only, since societal responsibility precludes such a luxury position. Also mathematics is facing increasingly complex problems, whose solution and subsequent solution-verification will require techniques beyond traditional practice. Examples (i) and (ii) above are just some first witnesses of this kind (in fact, there have been other examples before, including the four color theorem \cite{1977SciAm.237d.108A,Gonthier08}). Formal proofs therefore should, if not must, adopt 
a more central role, in mathematics and beyond.

However, vice versa, I clearly also argue for coupling formal proofs with additional human-intuitive proofs whenever feasible. Explainability, transparency and intuition must remain virtues of highest priority, not only in mathematics, but in particular in  topical, emerging areas such as autonomous intelligent machines. I am thus against preferring one notion over the other. Instead, both notions of proof should be coupled pari passu, whenever possible. And in the long run, the raised trustworthiness and beauty of a combined approach will justify the required additional resource efforts.

Note that I have avoided the 
phrase \emph{integrating mathematical proofs and formal proofs}. Instead, I have talked about \emph{pairings} or \emph{couplings} only. While I do not rule out that a proper integration of both notions of proofs can eventually be achieved (and there has been relevant research in the past, see e.g. \cite{DBLP:conf/tableaux/Bundy98,B10} and the references therein), 
considerable scientific progress is still needed to get there. The overall challenge, however, appears \emph{AI-complete}, since it includes a seamless, semantical integration of natural language, diagrams and formula language.

My own research is developing and applying pairings of formal proof and human-intuitive proof in an inspiring novel direction not mentioned so far: \emph{computational metaphysics}. In collaboration with colleagues, I have demonstrated that also in metaphysics (and ethics and argumentation) formal poofs have a lot to contribute, including the revelation of philosophically relevant new insights~\cite{C55}. For example, my higher-order theorem prover LEO-II \cite{J30} revealed an unnoticed inconsistency in G{\"o}del's modern variant of the ontological argument for the existence of God, while Scott's emendation of G{\"o}del's argument (and the consistency of the emended premises) was automatically verified. These applications in metaphysics were enabled by a new, generic technique in which \emph{classical higher-order logic (HOL)}, as supported in modern theorem provers and proof assistants, \emph{is utilized as a universal meta-logic} in which different target logics can be semantically embedded. In metaphysics we can thus encode arguments in \emph{higher-order modal logic} (as e.g.~assumed by G{\"o}del for his ontological argument), while in category theory we may want to work with \emph{free first-order logic}, which is suitable for addressing partiality issues in a proper way \cite{J40}. With the new technique, existing theorem provers for HOL become readily  applicable in all these application contexts \cite{J41}. 

% I am thus convinced that this technique is well suited to address a wide range of challenge applications for which other implementations of formal proofs are not available yet; see \cite{J41} for more on this.

So what is needed to further develop and foster the utilization of an integrated notion of formal and mathematical proof in future applications?
It is a next generation of highly qualified experts that master both, the beauty and intuition of mathematical proofs and the technicality challenges of formal proofs. Unfortunately, however, this vision has not been picked up yet by mathematicians to the extend that the deductions systems community was hoping for (see e.g.~the discussions in \cite{QED,Notices2008,E7,Bundy05}). 

The recent work by Marco David, Benjamin Stock, Abhik Pal and their fellow students at Jacobs University, however,
%as discussed in more detail in this article 
provides good new hope. 
While their ongoing verification project  \cite{EasyChair:152} on Matiyasevich's proof of  Hilbert's tenth problem is in many ways related to the Flyspeck project, albeit on a smaller scale, there is also a significant difference which I personally find particularly encouraging. While Hales, to my best knowledge, received substantial support already early on by expert members from the deduction system community, the maths students were entering their formal encoding project without prior knowledge of the proof assistant technology they employed (Isabelle/HOL \cite{Isabelle}). They initially also had little knowledge about the logical foundations of that system and for a long time into the project they worked without any expert support. And yet, they still mastered the challenge and got a very long way by working with the proof assistant on their own. This provides good evidence for the maturity proof assistant technology has meanwhile achieved. For a next generation of 
talented maths students, the required expertise acquisition can obviously be handled autonomously, while for the majority of more matured/established mathematicians this may not pose a realistic and sufficiently attractive scenario anymore.

To conclude, I am convinced that an integrated notion of formal and intuitive mathematical proof is indispensable for a wide range of topical, future applications across disciplines, and there is encouraging recent evidence that this is practically feasible.
Society cannot afford to postpone further research and developments in this area, with or without the support of traditional mathematics. Whether and when mathematicians are wholeheartedly embracing the new developments is a subordinate but nevertheless relevant question, since their support could trigger a substantial increase of much needed research and infrastructure investments.

\bibliographystyle{abbrv}
%\bibliography{bibliography}

\begin{biog}
\item[Christoph Benzm{\"u}ller] 
is a professor in computer science and mathematics at Freie Universit{\"a}t Berlin (Germany). He is also a permanent visiting scholar of the University of Luxembourg (Luxembourg). Christoph's prior research institutions include the universities of Stanford (USA), Cambridge, Birmingham, Edinburgh (all UK), the Saarland (Germany) and CMU (USA).

Christoph received his PhD (1999) and his Habilitation (2007) from Saarland University, his PhD studies were partly conducted at CMU. In 2012, Christoph was awarded with a Heisenberg Research Fellowship of the German National Research Foundation (DFG).

The research activities of Christoph are interfacing the areas of artificial intelligence, philosophy, mathematics, computer science, and natural language. Many of these activities draw on classical higher-order logic (HOL), which has is roots in the work by Russel and Church. Christoph has contributed to the semantics and proof theory of HOL, and together with colleagues and students he has developed the Leo theorem provers for HOL. More recently he has been utilizing HOL as a universal meta-logic to automate various non-classical logics in topical application areas, including machine ethics \& machine law (responsible AI), rational argumentation, metaphysics, category theory, etc. 
%Recent research also addresses the integration of automated reasoning and machine learning.
\begin{affil}
Department of Mathematics and Computer Science, Freie Universit{\"a}t Berlin, Arnimallee 7, 14195 Berlin\\
E-mail: c.benzmueller@fu-berlin.de
\end{affil}
\end{biog}
\vfill\eject

\end{document}